\input amstex
\documentstyle {amsppt}
\magnification=\magstep1
\hsize=16truecm
\vsize=22.5truecm
\baselineskip=16truept
\NoRunningHeads
\NoBlackBoxes
\loadbold

\topmatter
\title Desingularized fiber products of semi-stable elliptic surfaces with
vanishing third Betti number\endtitle
\author Chad Schoen$^1$ \endauthor
\abstract Desingularized fiber products of semi-stable elliptic surfaces with  $H^3_{etale}=0$ are classified.  Such varieties may play a role in the study of  supersingular threefolds, of the deformation theory of varieties with trivial canonical bundle, and of arithmetic degenerations of 
rigid Calabi-Yau threefolds.\endabstract
\address Department of Mathematics, Duke University,
Box 90320, Durham, NC 27708-0320 USA.\qquad
e-mail: schoen\@ math.duke.edu .\endaddress
\thanks $^1$Partial support by the NSF (DMS-0200012) and NSA (MDA904-97-1-0041) gratefully acknowledged \endthanks
\subjclass\nofrills{{\rm 2000}{\it Mathematics Subject Classification}.\usualspace } 14D06, 14J32, 14G15  \endsubjclass

\toc 
\head 1. Introduction \endhead
\head 2. Notations  \endhead
\head 3. The third Betti number of desingularized fiber products \endhead
\head 4. Semi-stable elliptic surfaces over $\Bbb P^1$ with four singular fibers \endhead
\head 5. Automorphisms of $\Bbb P^1$ which stabilize a set of four points \endhead
\head 6. Existence and non-existence of isogenies \endhead 
\head 7. Classification of desingularized fiber products with $h^3=0$ \endhead
\head 8. The canonical sheaf and some Hodge numbers \endhead
\head 9. Proof of Theorem 1.1 \endhead 
\head 10. Projective threefolds with trivial canonical sheaf and $h^3=0$ \endhead 
\head 11. Deformations and lifting to characteristic zero \endhead
\head 12. Supersingular threefolds \endhead 
\head 13. Arithmetic degeneration of rigid Calabi-Yau threefolds \endhead 
\head 14. Fiber products of more general elliptic surfaces \endhead 
\endtoc
\endtopmatter

\document

\def\hr{height2pt&\omit&&\omit&&\omit&&\omit&&\omit&&\omit&\cr}
\define\op{\operatorname}

\define\ql{\Bbb Q_l}
\define\inv{^{-1}}

\subhead 1. Introduction \endsubhead

Let $X$ be a smooth, projective, irreducible curve over an algebraically closed field $k$.
Let $\pi :Y\to X$ and $\pi ':Y'\to X$ be relatively minimal elliptic surfaces with section.
We assume that each is semi-stable and has at least one singular fiber. The fiber product,
$\bar W=Y\times _XY'$, is singular exactly at points $(y,y')\in \pi \inv (s)_{sing}\times (\pi ')\inv (s)_{sing}$. Blowing up the reduced singular locus of $\bar W$ 
yields a non-singular projective threefold, $W$. 
This note is concerned with such threefolds $W$ for which the third 
$l$-adic Betti number, $h^3(W,\ql )$, is zero ($l\neq \op{char}(k)$). 
That such threefolds exist at all may come as a surprise. 
The fact that the Hodge number, $h^0(W,\Omega ^3_{W/k})$, 
always turns out to be positive implies that they do not exist when 
$\op{char}(k)=0$.  A first step towards understanding these
objects is the following classification theorem.

\proclaim {Theorem 1.1} Let $W$ be constructed as above. Suppose that $h^3(W,\ql )=0$. Then

(i) $char(k)\in \{ 2,3,5,7,11,17,29,31,41,73,251,919,9001\} $.

(ii) $W$ may be defined over the prime subfield of $k$ unless $\op{char}(k)=2$,
in which case $W$ may be defined over the field with four elements.

(iii) For each characteristic $p$ in (i) and each pair $(n,n')\in (\Bbb Z_{\geq 0})^2$ there exists $W$ as above with $h^0(W,\Omega ^3_{W/k})=p^n+p^{n'}-1$. Other values of $h^0(W,\Omega ^3_{W/k})$ do not occur.

(iv) There are only finitely many isomorphism classes of
varieties $W$ as above with fixed Hodge number $h^0(W,\Omega ^3_{W/k})$.
\endproclaim

The proof of the theorem begins with the formula for the third \'etale Betti number of a desingularized fiber 
product of semi-stable elliptic surfaces with section given in 
Proposition 3.1. The Betti number is zero precisely when $X$ has genus $0$,
$\pi $ and $\pi '$ are not isogenous and have the same four places of bad 
reduction (cf. Corollary 3.2). The classification of semi-stable elliptic
surfaces over $\Bbb P^1$ with four singular fibers is recalled in \S 4.
The task of arranging that the places of bad reduction coincide without
forcing the elliptic fibrations to be isogenous is the focus of sections 
$5$ and $6$. The desingularized fiber products are classified in \S 7. 
Some Hodge numbers are computed in \S 8 and Theorem 1.1 is proved in \S 9.  
In section 10 resolutions of the singularities of $\bar W$ which do not
introduce exceptional divisors are discussed. This leads to the construction
of certain smooth projective threefolds, $\widehat W$, with trivial canonical sheaf and $H^3(\widehat W,\ql )=0$. Deformations of these varieties are
investigated in \S 11. Other smooth projective threefolds with trivial 
canonical sheaf and $H^3(\widehat W,\ql )=0$ were constructed by 
Hirokado \cite {Hi1}, Schr\"oer \cite {Schr}, and Ekedahl \cite {Ek} 
using very different methods. Some of the differences and the questions
that they raise are discussed in sections 12 and 13.

\example{Acknowledgement}
This project began when the author considered the reduction at various primes
of the varieties studied by M. Sch\"utt in \cite {Sch\"u}.
The author thanks Sch\"utt for correspondence and discussions in 2003 and 2004
which played an important role in this work. 
\endexample 

\bigpagebreak
\subhead 2. Notations \endsubhead 

\noindent $\Bbb N=\Bbb Z_{\geq 1}$.

\noindent $k=$ an algebraically closed field.

\noindent $l=$ a prime number different from $\op{char}(k)$.

\noindent $X=$ an irreducible curve, smooth and projective over $k$.

\noindent $g_X=$ the genus of $X$.

\noindent $j:\eta =Spec(k(X))\to X$ the inclusion of the generic point of $X$.

\noindent $\bar {\eta }=Spec(\overline {k(X)})$, a geometric generic point which is algebraic over $\eta $.

\noindent $G_{\eta }=Gal(\bar {\eta }/\eta )$ the absolute Galois group of the field $k(X)$.

\noindent Semi-stable elliptic surfaces are assumed to be relatively minimal and to have 
a section. 

\noindent $\pi :Y\to X$, $\pi ':Y\to X$ non-isotrivial semi-stable elliptic surfaces.

\noindent $\bar W=Y\times _XY'$.

\noindent $W@>\sigma >>\bar W$ is the blow-up of $\bar W$ along the reduced singular locus.

\noindent $Q\subset W$ is the exceptional divisor for $\sigma $.

\noindent $\bar f:\bar W\to X$ is the tautological map. $f:=\bar f\circ \sigma :W\to X$. 

\noindent $S$ (respectively $S'$) is the locus in $X$ over which $\pi $ (respectively $\pi '$)
fails to be smooth.

\noindent $m_s$ (respectively $m'_s$) number of irreducible components of $\pi \inv (s)$
(respectively $(\pi ')\inv (s)$).

\noindent Two elliptic surfaces over a common base curve are {\it isogenous} if 
their generic fibers are.

\noindent $\epsilon =\cases 1,& \ \ \text{if} \ Y \ \text{and} \ Y' \ \text{are isogenous} \\ 0,& \ \ \text{otherwise}. \endcases $

\bigpagebreak
\subhead 3. The third Betti number of desingularized fiber products      \endsubhead 

With notation defined in the previous section:

\proclaim {Proposition 3.1} $h^3(W,\ql )=2(6g_X-4+|S\cap S'|+\epsilon +\sum _{s\in S-S\cap S'}m_s+\sum _{s\in S'-S\cap S'}m'_s)$.
\endproclaim 

\proclaim {Corollary 3.2} $h^3(W,\ql )=0$ exactly when $g_X=0$, $\epsilon =0$, 
$S=S'$ and $|S|=4$.
\endproclaim
\demo{Proof} By the proposition $h^3(W,\ql )=0$ implies $g_X=0$.
According to \cite {Be1} the minimal number of singular fibers of 
a semi-stable elliptic surface over $\Bbb P^1$ is $4$.  Thus $|S\cup S'|=|S\cap S'|+|S-S\cap S'|+|S'-S\cap S'|\geq 4$ with equality only when $S=S'$ and $|S|=4$.
The assertion follows.
\hfill $\square $\enddemo 

\demo{Proof of Proposition 3.1}
The Betti number $h^3(W,\ql )$ is independent of the choice of prime $l\neq \op{char}(k)$
\cite {Mi, VI.12.5b}. The proof proceeds via the Leray spectral sequence for $f:W\to X$. Applying Poincar\'e duality \cite {Mi, V.2.2c} and ignoring Tate twists (since $k$ is algebraically closed) gives
$$E_2^{2,1}\simeq H^2(X,R^1f_*\ql )\simeq H^2(X,j_*j^*R^1f_*\ql )\simeq H^0(X,j_*j^*R^1f_*\ql ^{\vee })=0,$$
since 
$$H^1(f\inv (\bar {\eta }),\ql )^{\vee }\simeq H^1(\pi \inv (\bar {\eta }),\ql )\oplus H^1((\pi ') \inv (\bar {\eta }),\ql )$$
has no non-zero $G_{\eta }$-invariants. Indeed, the existence of non-zero $G_{\eta }$-invariants would imply the existence of a dominant morphism from $X$ to a tower of modular curves, $X_1(l^n)$, $n\to \infty $, of unbounded genus which is impossible.

Observe that  $E_{\infty }^{1,2}$ is isomorphic to 
$$E_2^{1,2}\simeq H^1(X,R^2f_*\ql )\simeq H^1(X,j_*j^*R^2f_*\ql ),$$
where the second isomorphism is a consequence of the local invariant cycle theorem
\cite {De,3.6.1}. Set $\frak R:=j_*j^*(R^1\pi _*\ql \otimes 
R^1\pi '_*\ql )$ and note that $j_*j^*(R^2\pi _*\ql \otimes \pi '_*\ql )\simeq \ql (-1) .$
Thus 
$$j_*j^*R^2f_*\ql \simeq \ql (-1) \oplus \ql (-1) \oplus \frak R.\tag 3.1$$
Since $\frak R$ is self-dual up to a Tate twist, Poincar\'e duality \cite {Mi, V.2.2c} yields,
$$h^1(X,\frak R)=-e(X,\frak R)+2h^0(X,\frak R).$$
The isogeny theorem gives $h^0(X,\frak R)=\epsilon $ 
\cite{Fa-Scha-W\"u, IV.1 Remark (ii)}. There is a 
formula for the euler characteristic,
$$e(X,\frak R)=e(X,\ql )dim_{\ql }(\frak R_{\bar {\eta }}) -
\sum _{s\in (S\cup S')}\left[ dim_{\ql }(\frak R_{\bar {\eta }})-dim_{\ql }(\frak R_{\bar{\eta }}^{I_s})\right] ,\tag 3.2$$
which will be verified in Lemma 3.3 below for certain primes $l$.
Now $dim_{\ql }(\frak R_{\bar {\eta }})=4$.
It follows easily from the known action of the inertia group at a fiber of 
multiplicative reduction \cite {Sil2,V.Ex5.13} that 
$dim_{\ql }((\frak R_{\bar{\eta}})^{I_s})=2$
for any $s\in S\cup S'$. Thus (3.2) gives 
$$h^1(X,\frak R)=-e(X,\frak R)+2h^0(X,\frak R)=-(4(2-2g_X)-2|S\cup S'|)+2\epsilon .$$
Adding in the two $h^1(X,\ql (-1))$ terms from (3.1) yields
$$h^1(X,j_*j^*R^2f_*\ql )=4g_X-(4(2-2g_X)-2|S\cup S'|)+2\epsilon .\tag 3.3$$

Consider finally the term, $E_2^{0,3}= H^0(X,R^3f_*\ql )$. By the local invariant 
cycle theorem the map on stalks, 
$$(R^3f_*\ql )_s\to (j_*j^*R^3f_*\ql )_s\simeq (j_*j^*(R^1\pi _*\ql \otimes R^2\pi '_*\ql \ \oplus \ R^2\pi _*\ql \otimes R^1\pi '_*\ql ))_s ,\tag 3.4$$
is surjective for any closed point $s\in X$. For $s\in S\cap S'$ 
$h^3(f\inv (s),\ql )=2$  \cite {Sch3,\S 12} and 
the target has dimension $2$, since $dim_{\ql }(R^1\pi _*\ql )_{\bar{\eta }}^{I_s}=1$ at a place of multiplicative reduction and the same holds for $\pi '$. 
The proper base change theorem implies that (3.4) is an isomorphism. As $H^0(X,j_*j^*R^3f_*\ql )=0$, the global sections of $R^3f_*\ql $ are all supported on $(S\cup S')-(S\cap S')$. 
Set $U=X-(S\cap S')$ and observe that  
$$R^3f_*\ql |_U\simeq R^3\bar f_*\ql |_U\simeq (R^1\pi _*\ql \otimes R^2\pi '_*\ql \ \oplus \ R^2\pi _*\ql \otimes R^1\pi '_*\ql )|_U$$
decomposes as 
$$(R^1\pi _*\ql \otimes (\ql (-1) \oplus \sum _{s\in S'\cap U}i_{s*}\ql ^{m_s'-1}))|_U\oplus ((\ql (-1)\oplus \sum _{s\in S\cap U}i_{s*}\ql ^{m_s-1})\otimes R^1\pi '_*\ql )|_U.$$
Thus 
$$h^0(X,R^3f_*\ql )=\sum _{s\in S-(S\cap S')}2(m_s-1)\ + \ 
\sum _{s\in S'-(S\cap S')}2(m'_s-1).\tag 3.5 $$
Furthermore the differential $d_2^{0,3}:E_2^{0,3}\to E_2^{2,2}$ is zero.
This may be seen by noting that $f$ is defined over some finitely generated
$\Bbb Z$-algebra, specializing to a suitable place with finite residue field
and applying a weight argument \cite {Ka, 7.5.2}.

Thus $h^3(W,\ql )=h^1(X,j_*j^*R^2f_*\ql ) + h^0(X,R^3f_*\ql )$ and the proposition
follows from equations (3.3) and (3.5).
\hfill $\square $\enddemo
It remains to justify equation (3.2).
\proclaim {Lemma 3.3} For primes $l\neq \op{char}(k)$ and $l\nmid 2\prod _{s\in S}m_s\prod _{s\in S'}m'_s$, 
$$e(X,\frak R)=e(X,\ql )dim_{\ql }(\frak R_{\bar {\eta }}) -
\sum _{s\in (S\cup S')}\left[ dim_{\ql }(\frak R_{\bar {\eta }})-dim_{\ql }(\frak R_{\bar{\eta }}^{I_s})\right] .$$
\endproclaim
\demo{Proof} For an abelian group $A$, let $A[l]$ denote the kernel of multiplication
by $l$. Define $\frak R_n:= j_*j^*(R^1\pi _*\Bbb Z/l^n\otimes R^1\pi '_*\Bbb Z/l^n)$ and
set $\varpi ^i:=dim_{\Bbb Z/l}\left( (\varprojlim _n H^i(X,\frak R_n)) [l]\right) $.
Since $\pi $ and $\pi '$ are semi-stable the action of any inertia group
on the stalk $\frak R_{1,\bar{\eta}}$ factors through the tame quotient \cite {Sil2, IV.10.2b}
and \cite {Mi, V.2.12} gives the Euler characteristic formula,
$$e(X,\frak R_1)=e(X,\Bbb Z/l)dim_{\Bbb Z/l}(\frak R_{1\ \bar {\eta }}) -
\sum _{s\in (S\cup S')}\left[ dim_{\Bbb Z/l }(\frak R_{1\ \bar {\eta }})-dim_{\Bbb Z/l }(\frak R_{1\ \bar{\eta }}^{I_s})\right] .\tag 3.6$$
The hypothesis on $l$ implies that the right hand sides of (3.6) and (3.2) are equal. 
and that each $\frak R_n$ is a flat sheaf of $\Bbb Z/l^n$-modules.
From the exact sequence \cite {Mi, proof of V.1.11},
$$\varprojlim _nH^i(X,\frak R_n))@>l>>\varprojlim _nH^i(X,\frak R_n))
@>>>H^i(X,\frak R_1)@>>>\varprojlim _nH^{i+1}(X,\frak R_n))@>l>>
\varprojlim _nH^{i+1}(X,\frak R_n)),$$
one deduces a short exact sequence,
$$0@>>>(\Bbb Z/l)^{h^i(X,\frak R)+\varpi ^i}@>>>H^i(X,\frak R_1)@>>>(\Bbb Z/l)^{\varpi ^{i+1}}
@>>>0,$$
for each $i$. The equality $e(X,\frak R_1)=e(X,\frak R)$ and with it the lemma follow. 
\hfill $\square $\enddemo

\bigpagebreak
\subhead 4. Semi-stable elliptic surfaces over $\Bbb P^1$ with four singular fibers \endsubhead

Semi-stable elliptic surfaces over $\Bbb P^1_k$ with four singular fibers 
were classified by Beauville \cite {Be2} when $\op{char}(k)=0$. They are all 
rational surfaces coming from pencils of cubic 
curves on $\Bbb P^2$. The table below contains Beauville's description of the pencil,
the set $S\subset \Bbb P^1$ where the pencil has bad reduction and the Kodaira type 
of the singular fiber above each point of $S$. Each fibration is specified by 
an integer which we call the {\it level} and denote by $L$ in the table. An identity
section defined over $\Bbb Z[1/L]$ is specified in the third column by giving a 
base point for the pencil. This is a simple base point when $L\in \{ 3,5,6\} $.
When $L\in \{ 4,8\} $ the section comes from the unique infinitely near base point
lying over the given base point and when $L=9$ it comes from the second order
infinitely near base point. The next to
last column gives the cross ratio orbit (c.r.o.)  of $S$, that is the set of all $\lambda \in k$ such that there is an automorphism of $\Bbb P^1_k$ which
takes the set $S$ to $\{ 0,1,\infty ,\lambda \} $. If $\lambda $ is one such element,
then the others are  
$1/\lambda , 1-\lambda , 1/(1-\lambda ), \lambda /(\lambda -1), (\lambda -1)/\lambda $. 
The final column gives the $J$-invariants of the elliptic fibrations (cf. \cite {Sil, III.1}). These may be computed starting from the Weierstrass equations in \cite {Mi-Pe, Table 5.3}. 
\medskip
\proclaim{\hfil Table 4.1 Semi-stable rational elliptic surfaces with four singular fibers}
$$\eightpoint{        
\vbox {\offinterlineskip\tabskip=0pt\hrule
\halign{ &\vrule#&\strut $\quad #\hfil\quad$\cr
\hr
&\text{L} && \text{Equation} && \text{section} &&\text{Fibers} && S && 
\text{c.r.o.} && J&\cr
\hr
\noalign{\hrule}
\hr
& 3 && \ssize{x^3+y^3+z^3+txyz} && \ssize{(1:-1:0)} && \text{I}_3 \text{I}_3\text{I}_3\text{I}_3&& \ssize{-3\mu _3, \infty } && \zeta _6, \zeta _6\inv  && \frac{(t^4-6^3t)^3}{(t^3+27)^3}&\cr
\hr
\noalign{\hrule}
\hr
& 4 && \ssize{x(x^2+y^2+2zy)+tz(x^2-y^2)} && \ssize{(0:0:1)} && \text{I}_2 \text{I}_2\text{I}_4\text{I}_4&& \ssize{1,-1,0,\infty }&& -1,2,\frac 12 && \frac{-2^8(t^4-t^2+1)^3}{t^4(t+1)^2(t-1)^2}&\cr
\hr
\noalign{\hrule}
\hr
& 5 && \ssize{x(x-z)(y-z)+tzy(x-y)} && \ssize{(1:0:1)} && \text{I}_1 \text{I}_1\text{I}_5\text{I}_5&& \ssize{\Theta ,0,\infty }&& \frak K && \frac{(t^4-12t^3+14t^2+12t+1)^3}{t^5(t^2-11t-1)}&\cr
\hr
\noalign{\hrule}
\hr
& 6 && \ssize{(x+y)(y+z)(z+x)+txyz} && \ssize{(1:-1:0)}&& \text{I}_1 \text{I}_2\text{I}_3\text{I}_6&& \ssize {-8, 1, 0, \infty} && \frak C &&\frac{-(t^4+8t^3-16t+16)^3}{t^3(t-1)^2(t+8)}&\cr 
\hr
\noalign{\hrule}
\hr
& 8 && \ssize{(x+y)(xy-z^2)+txyz} && \ssize{(0:0:1)}&& \text{I}_1 \text{I}_1\text{I}_2\text{I}_8&& \ssize{-1, 1, 0, \infty} && -1,2,\frac 12&& \frac{2^4(16t^4-16t^2+1)^3}{t^2(t+1)(t-1)}&\cr 
\hr
\noalign{\hrule}
\hr
& 9 && \ssize{x^2y +y^2z+z^2x+txyz} && \ssize{(0:0:1)} && \text{I}_1 \text{I}_1\text{I}_1\text{I}_9&& \ssize{-3\mu _3, \infty }&& \zeta _6, \zeta _6\inv && \frac{t^3(t^3+24)^3}{t^3+27}&\cr 
\hr
\noalign{\hrule}
\hr
}\hrule}                } $$
\endproclaim
In the fifth column $\mu _3$ denotes the set of all third roots of one and 
$\Theta =\{ \theta ,\theta '\} $ indicates the roots of the polynomial $t^2-11t-1$. In the sixth 
column $\zeta _6$ is a primitive cube root 
of $-1$, $\frak C=\{ -8,-\frac 18, 9, \frac 19, \frac 98, \frac 89\} $ and 
$\frak K=\{ \kappa,1/\kappa ,1-\kappa , 1/(1-\kappa ), (\kappa -1)/\kappa ,
\kappa /(\kappa -1)\} $ where $\kappa =\theta /\theta '$ ($\kappa =(11+5\sqrt 5)/(11-5\sqrt 5)$ if $char(k)\neq 2$).

\proclaim{Proposition 4.1} The above table gives a complete classification of rational, semi-stable  elliptic surfaces with four singular fibers in any characteristic 
with the one caveat that no surface of level $L$ exists in characteristic $p$ when $p|L$. 
\endproclaim
\demo{Proof} W. Lang \cite {La}.
\hfill $\square $\enddemo
Non-rational semi-stable elliptic surfaces over $\Bbb P^1_k$  with four 
singular fibers may be constructed in  characteristic $p>0$ by Frobenius 
base change. We recall the construction: Each entry in Table 4.1 with $p\nmid L$
determines an elliptic surface over the prime field, $\pi _0:Y_o\to \Bbb P^1_{\Bbb F_p}$. 
Let $F_0:\Bbb P^1_{\Bbb F_p}\to \Bbb P^1_{\Bbb F_p}$ be the absolute Frobenius
morphism, ie. the identity on the topological space and the $p$-th power 
map on the sturcture sheaf. For each $n\in \Bbb Z_{\geq 0}$ there is a commutative 
diagram,
$$\CD 
Y_{0,n}@>\varrho _{0,n}>>\bar Y_{0,n} @>>> Y_0 \\
@V\pi _{n,0}VV @V\bar {\pi }_{n,0}VV @V\pi _0VV \\
\Bbb P^1_{\Bbb F_p}@=\Bbb P^1_{\Bbb F_p}@>F_0^n>>P^1_{\Bbb F_p}.
\endCD $$
in which the right hand square is Cartesian and $\varrho _{0,n}$ is a 
minimal resolution of singularities. Base change the diagram 
by $Spec(k)\to Spec(\Bbb F_p)$ and write $X=\Bbb P^1_k$:
$$\CD 
Y_n@>\varrho _n>>\bar Y_n @>>> Y \\
@V\pi _nVV @V\bar {\pi }_nVV @V\pi VV \\
X@=X@>F^n>>X.
\endCD \tag 4.1$$
\proclaim {Lemma 4.2}
(i) $\pi _n:Y_n\to X$ is an elliptic surface with $4$ singular fibers.

(ii) The critical loci of $\pi _n$ and $\pi $ coincide. 

(iii) To each singular fiber of $\pi $ of type I$_r$ there is a corresponding
singular fiber of $\pi _n$ of type I$_{p^nr}$.

(iv) The elliptic fibrations $\pi _n$ and $\pi $ are isogenous.
\endproclaim
\demo{Proof} (i) Let $S\subset X$ be the complement of the locus over
which $\pi $ is smooth. Then $|S|=4$ and $\pi _n$ is smooth over $X-(F^n)\inv (S)$.  

(ii) Since $S$ is defined over the prime field, 
$S=(F^n)\inv (S)$.

(iii) The corresponding fiber of $\bar {\pi }_n$ is an $r$-gon of $\Bbb P^1$'s.
A local calculation gives that the surface $\bar Y_n$ has an $A_{p^n-1}$
singularity at each singular point of the $r$-gon. Blowing these up gives 
the result.

(iv) The $n$-th power of the absolute Frobenius map, $F_{Y_0}^n:Y_0\to Y_0$, 
factors through the fiber product, $\bar Y_{0,n}$. The 
identity section maps to the identity section so the map on generic
fibers is an isogeny. Now base changing gives a morphism, 
$Y\to \bar Y_n$, which is an isogeny on generic fibers.
\hfill $\square $\enddemo 

\example {Definition 4.3} Fix an algebraically closed field $k$. The elliptic 
fibrations which appear in Table 4.1 will be refered to as BLIS fibrations
of exponent $0$. In the event that $char(k)>0$, the elliptic fibration, $\pi _n$, 
constructed in (4.1) from the BLIS fibration of exponent $0$, $\pi $, 
will be called a BLIS fibration of exponent $n$.
\endexample 

\proclaim {Proposition 4.4}
Every semi-stable elliptic fibration over $\Bbb P^1_k$ with exactly $4$ singular 
fibers is isomorphic to a unique BLIS fibration.
\endproclaim
\demo{Proof} 
Semi-stable elliptic surfaces over $\Bbb P^1_k$ with four singular fibers are {\it extremal}
in the sense that the Picard number equals the second Betti number and the Mordell-Weil 
group is finite \cite {Be1, Section 4.A}. Extremal elliptic surfaces were classified by
H. Ito \cite {I1}, \cite {I2}, and A. Schweizer \cite {Schw}.
\hfill $\square $\enddemo

\example {Remark 4.5} BLIS is short for Beauville-Lang-Ito-Schweizer.
\endexample

\example {Remark 4.6} The elliptic surfaces of levels $3$ and $9$ in Table 4.1
are isogenous. The same holds for the surfaces of levels $4$ and $8$.
\endexample

\bigpagebreak
\subhead 5. Automorphisms of $\Bbb P^1$ which stabilize a set of four points \endsubhead

It is important to study automorphisms of $\Bbb P^1$ which stabilize a given set of four points
for the following reason: Suppose that $\pi :Y\to X=\Bbb P^1_k$ is a BLIS fibration,  
$\alpha \in Aut(X)$ stabilizes the bad reduction locus $S$, and the pullback, $\pi '$
of $\pi $ by $\alpha $ is not isogenous to $\pi $. Then the desingularized fiber
product of $\pi $ and $\pi '$ has $h^3=0$ by Corollary 3.2.

Given a finite subset of closed points, $T\subset \Bbb P^1_k$, $Aut(T)$ will denote 
the set of permutations of the set while $Aut(\Bbb P^1_k,T)$ denotes the subgroup of 
$Aut_k(\Bbb P^1_k)$ which stablizes $T$. If $|T|\geq 3$, then $Aut(\Bbb P^1_k,T)$
is canonically identified with a subgroup of $Aut(T)$. Assume henceforth that
$|T|=4$. Write $V_4\subset A_4\subset \frak S_4$ for Klein's fourgroup, the 
alternating group, respectively the symmetric group. Let $D_4$ denote the 
dihedral group with $8$ elements.
When $char(k)\neq 3$ write $\zeta _6\in k^*$ for a primitive cube root of $-1$.

\proclaim {Lemma 5.1}
(i) If $char(k)\neq 2$ and the cross ratio orbit of $T$ is $\{ -1,2,1/2\} $, then 
$Aut(\Bbb P^1_k,T)\simeq D_4$, except when $char(k)=3$, in which case 
$Aut(\Bbb P^1_k,T)\simeq \frak S_4$. 

(ii) If $char(k)\neq 3$ and the cross ratio orbit of $T$ is $\{ \zeta _6, \zeta _6\inv \} $, then $Aut(\Bbb P^1_k,T)\simeq A_4$.

(iii) If the cross ratio orbit of $T$ is different from cases (i) and (ii), then
it consists of $6$ distinct points and $Aut(\Bbb P^1_k,T)\simeq V_4$.
\endproclaim
\demo{Proof} For $\lambda \in k^*-\{ 1\} $, the subgroup, $\frak V\subset Aut_k(\Bbb P^1_k)$, generated by $t\mapsto \lambda /t$ and $t\mapsto (t-\lambda )/(t-1)$ 
is isomorphic to $V_4$ and stabilizes $\{ \infty ,0,1,\lambda \} $. 
Any element of $Aut_k(\Bbb P^1_k)$ which maps $\{ \infty ,0,1,\lambda \} $
to a set containing $\{ \infty ,0,1\} $ has the form $\frak s\circ \frak v$
where $\frak v\in \frak V$ and $\frak s\in Aut(\Bbb P^1_k, \{ \infty , 0,1\} )
\simeq \frak S_3$. If only the identity in $Aut(\Bbb P^1_k, \{ \infty , 0,1\} )$
fixes $\lambda $, then the cross ratio orbit consists of $6$ distinct points
and $Aut(\Bbb P^1_k,\{ \infty ,0,1,\lambda \} )\simeq V_4$. The conjugacy class in 
$Aut(\Bbb P^1_k, \{ \infty , 0,1\} )$ of points of order $2$ consists of 
the fractional linear transformations $t\mapsto 1/t$ (respectively $t/(t-1)$, $1-t$) 
which fix the pairs of points $-1,1$ (respectively $2,0$;  $1/2,\infty $). 
Case (i) when $char(k)\neq 3$ follows. The conjugacy class of elements
of order $3$ consists of the transformations, $t\mapsto 1/(1-t)$ 
(respectively $(t-1)/t$). These automorphisms 
fix only $\zeta _6$ and $\zeta _6\inv $ and (ii) follows. Finally
when $char(k)=3$, $GL(\Bbb F_3)/\pm Id\simeq \frak S_4$ acts faithfully on 
$\{ \infty ,0,1,-1\} $.
\hfill $\square $\enddemo

\example {Definition 5.2} The {\it type} of a semi-stable elliptic 
surface over $X=\Bbb P^1_k$ with four singular fibers is the product of the 
distinct primes which divide the level of the corresponding BLIS fibration
of exponent $0$.
\endexample

\example {Remark 5.3} (i) The type is an element of the set $\{ 2,3,5,6\} $. 

(ii) There are no elliptic surfaces of type $\tau $ in characteristic $p$ if $p|\tau $.

(iii) Two BLIS fibrations of the same type over algebraically closed fields of the
same characteristic $p\geq 0$  have the same set $S\subset \Bbb P^1(\Bbb F_{p^2})$ of bad reduction.
\endexample

\example {Definition 5.4} Write $Aut(\tau ,p)$ for $Aut(X,S)$, where $\pi :Y\to X=\Bbb P^1_k$
is any BLIS fibration of type $\tau $ over an algebraically closed field $k$ of 
characteristic $p$.
\endexample 

The canonical specialization map $\phi :Aut(\tau ,0)\to Aut(\tau ,p)$ is 
defined and injective for all primes $p\nmid \tau $ since $S$ in Table 4.1
is finite and \'etale over $Spec(\Bbb Z[1/\tau ])$. In fact $\phi $ is 
an isomorphism for all except finitely many characteristics $p$. The exceptional
characteristics are indicated in the following table. The notation 
$(G;p_1,...,p_r)$ means that $Aut(\tau ,p)\simeq G$, when $p\in \{ p_1,...,p_r\} $.

\medskip
\proclaim{\hfil Table 5.1 Automorphism groups}
$$\vbox
{\offinterlineskip\tabskip=0pt\hrule
\halign{ &\vrule#&\strut $\quad #\hfil\quad$\cr
\hr
&\tau  && Aut(\tau ,0) && Exceptions &\cr
\hr
\noalign{\hrule}
\hr
& 2 && D_4 && (\frak S_4;3) &\cr
\hr
\noalign{\hrule}
\hr
& 3 && A_4 && \text{none} &\cr
\hr
\noalign{\hrule}
\hr
& 5 && V_4&& (A_4;2, 31), (D_4;11,251) &\cr
\hr
\noalign{\hrule}
\hr
& 6 && V_4 && (D_4;5,7,17), (A_4;73) &\cr 
\hr
\noalign{\hrule}
\hr
}\hrule}$$
\endproclaim

For each exceptional value of $(\tau ,p)$ in Table 5.1 the next table gives 
automorphisms, $\alpha $, which represent the non-trivial left cosets of 
the image of $Aut(\tau ,0)\to Aut(\tau ,p)$ in $Aut(\tau ,p)$. The notation 
$\zeta _3$ refers to a primitive cube root of $1$; $S$ is as in Table 4.1.

\medskip
\proclaim{\hfil Table 5.2 Exceptional automorphisms}
$$\vbox
{\offinterlineskip\tabskip=0pt\hrule
\halign{ &\vrule#&\strut $\quad #\hfil\quad$\cr
\hr
&\text{characteristic}  && \tau  && S && \alpha &\cr
\hr
\noalign{\hrule}
\hr
& 2 && 5 && \infty , \mu _3 && t\mapsto \zeta _3t; t\mapsto \zeta _3^2t &\cr
\hr
\noalign{\hrule}
\hr
& 3 && 2 && \infty ,0,1,-1 && t\mapsto 1/(1-t); t\mapsto (t-1)/t &\cr 
\hr
\noalign{\hrule}
\hr
& 5 && 6 && \infty ,0,1,-8 && t\mapsto t/(t-1) &\cr 
\hr
\noalign{\hrule}
\hr
&7 &&  6 && \infty ,0,1,-8 && t\mapsto 1/t &\cr 
\hr
\noalign{\hrule}
\hr
&11 && 5  && \infty ,0,1,-1 && t\mapsto 1/t &\cr 
\hr
\noalign{\hrule}
\hr
&17 &&  6 && \infty ,0,1,-8 && t\mapsto 1-t &\cr 
\hr
\noalign{\hrule}
\hr
&31&& 5  && \infty ,0,5,6 && t\mapsto (5t-25)/t; t\mapsto 25/(5-t) &\cr 
\hr
\noalign{\hrule}
\hr
&73&&  6 && \infty ,0,1,-8 && t\mapsto 1/(1-t); t\mapsto (t-1)/t  &\cr
\hr
\noalign{\hrule}
\hr
&251 && 5 && \infty ,0, 171, 91 && t\mapsto 91-t &\cr
\hr
\noalign{\hrule}
\hr
}\hrule}$$
\endproclaim

\example {Justification of Table 5.1}
The justification is computational, but easy. 
We consider briefly the case $\tau =5$ to give the flavor.
From Table 4.1 $S=\{ \infty ,0,\theta ,\theta '\} $, where 
$\theta $ and $\theta '$ are the roots of 
$t^2-11t-1$. Now $Aut(X,S)$ is strictly larger than $V_4$ precisely when 
$\lambda =\theta /\theta '$ lies in one of the two exceptional cross ratio
orbits, $\{ -1,2,1/2\} $ or $\{ \zeta _6, \zeta _6\inv \} $. Suppose 
$\lambda =\zeta _6$. Then $\theta ^3=-(\theta ')^3$ which implies 
$0=11\theta ^2+\theta +11(\theta ')^2 +\theta '$. As 
$\theta +\theta '=11$ and $\theta \theta '=-1$, the previous 
equation may be rewritten as $0=11(11^2+2)+11=2^2\cdot 11\cdot 31$.
Characteristic $11$ 
is ruled out, since then $\theta '=-\theta $ so $\lambda =-1$. 
Characteristics $2$ and $31$ really are exceptional as the explicit
automorphisms $\alpha $ in Table 5.2 show. The case 
$\lambda \in \{ -1,2,1/2\} $ is left to the reader.
\endexample 

\example {Remark 5.5} A modular interpretation of the exceptional 
automorphisms is not apparent. At least in characteristics $2$, $3$ 
and $5$ these automorphisms do not stabilize the supersingular loci ($t=1$, 
respectively $t=\pm \sqrt{-1}$, respectively $t\in \{ 3,4,\pm \sqrt 2\} $).
\endexample 

\bigpagebreak
\subhead 6. Existence and non-existence of isogenies \endsubhead

Let $\pi :Y\to X=\Bbb P^1_k$ be a BLIS fibration of type $\tau $. Set $p=char(k)$. 
An element $\alpha \in Aut(\tau ,p)$ is said to be {\it generic} if it is 
in the image of the specialization map, $\phi :Aut(\tau ,0)\to Aut (\tau ,p)$, otherwise
it is said to be {\it exceptional}. Write $\pi ':Y'\to X$ for the base change of $\pi $ 
with respect to the automorphism, $\alpha  :X\to X$.

\proclaim {Proposition 6.1} (i) If $\alpha $ is generic, then $Y'$ is isogenous
to $Y$.

(ii) If $\alpha $ is exceptional, then $Y'$ is not isogenous to $Y$.
\endproclaim
\demo{Proof} 
(i) By Lemmma 4.2(iv) it suffices to verify the assertion when $\pi $ 
is one of the fibrations in Table 4.1. 
By Remark 4.6 it suffices to treat the cases of levels $3$, $4$, $5$ and $6$.
Define $\dot X=X-S$. Adding a $\dot {}$ to the notation indicates base 
change by the inclusion, $\dot X\subset X$. For each entry in Table 4.1
$\dot {\pi }:\dot Y\to \dot X$ is the universal family of elliptic curves
with a particular level structure. In the following we will use standard 
notation from the theory of moduli of elliptic curves \cite {Sil, Appendix C \S 13}.

Level $3$. In this case $\dot {\pi }:\dot Y\to \dot X$ is the universal family of elliptic curves with a symplectic level $3$ structure.
$SL(2,\Bbb Z/3)$ permutes the symplectic level three structures transitively.
This action on the functor gives an action on the universal family,
$\dot {\pi }$. The action on the modular curve, $\dot X$, is via $SL(2,\Bbb Z/3)/\pm Id
\simeq A_4$. Thus in level $3$ the base changed fibration, $\pi '$, is 
actually isomorphic to $\pi $.

Level $4$. The relevant level structure is a point of order $4$ plus a 
point of order $2$ which is not a multiple of the point of order $4$. 
The automorphism group of $\Bbb Z/2\times \Bbb Z/4$ is isomorphic to 
a dihedral group, $D_4$. This group acts on the modular curve, $\dot X$, through
the quotient by its center, $D_4/\pm Id\simeq Gal (\dot X/X_1(4))\times Gal(\dot X/X(2))$.
The generator of the first factor fixes the places of type $I_2$ reduction and 
interchanges the places of type $I_4$ reduction. 
The generator of the second factor fixes the places of type $I_4$ reduction, but 
interchanges the places of type $I_2$ reduction. 
Taking the quotient of $Y$ by the subgroup generated by $2s$, where $s$ is any section
of order $4$, gives an isogeny from $Y$ to the pull back of $Y$ via an involution
which takes places of type $I_4$ reduction to places of type $I_2$ reduction.

Level $5$. The relevant level structure is a point of order $5$. The operation of 
$(\Bbb Z/5)^*$ permutes the possible level structures and acts on the the modular
curve $\dot X$ via $(\Bbb Z/5)^*/\pm Id \simeq Gal (\dot X/X_0(5))$. This action
interchanges the two places of type $I_5$ reduction as well as the two places of 
type $I_1$ reduction. Taking the quotient by the distinguished subgroup of order 
five gives rise to a model of the base change of $Y$ via an Atkin-Lehner involution which 
interchanges places of type $I_5$ reduction with places of type $I_1$ reduction.

Level $6$. The relevant level structure is a point of order $6$. There are 
four subgroups of the distinguished cyclic subgroup of sections of order $6$
of $\pi $. Taking the quotient of $Y$ by each of these yields a model for the 
base change of $\pi $ with respect to each of the elements of $V_4\simeq Aut(6,0)$.

\medpagebreak
(ii) Suppose now that $\alpha \in Aut(X,S)$ is an exceptional automorphism.
To show that $Y$ and $Y'$ are not isogenous it suffices by Lemma 4.2(iv) and 
Remark 4.6 to treat the case that $Y$ appears in Table 4.1 and has level 
$3$, $4$, $5$, or $6$. For type $3$ there are no exceptional automorphisms
and hence nothing to prove. In the remaining cases we claim that $\pi $
and $\pi '$ are not isomorphic. In all cases except the case $\tau =5$ and  
$\op{char}(k)=11$ this is clear from the Kodaira types of the fibers at the 
places of bad reduction (cf. Tables 4.1 and 5.2). The exceptional 
case is settled by a look at the $J$-invariant (again see Tables 4.1 and 5.2).
Thus if there is an isogeny, $\psi :Y\to Y'$, then there is an isogeny with non-trivial
cyclic kernel. We consider the various possibiities:

Suppose that $\psi $ exists and that $\op{Ker}(\psi )$ has a non-trivial 
connected component. Taking the 
quotient by this component gives a BLIS fibration of type $\pi _n$
with $n>0$ which has fibers of Kodaria type I$_{p^nr}$ (cf. Lemma 4.2).
Taking a further quotient by an \'etale group scheme of order prime to $p$
yields an elliptic surface with type I$_{p^nr'}$ reduction for some $r'$ prime to $p$.
Since $p$ is prime to the type, this surface cannot be isomorphic to $Y'$.
Contradiction.

It remains to rule out the case that $\op{Ker}(\psi )$ is a cyclic \'etale
group scheme. Certain cyclic groups of this type occur naturally as part
of the level structure in the elliptic surfaces of Table 4.1. By the 
proof of part (i) taking the quotient of $Y$ by such a 
subgroup yields the pullback of $Y$ by a generic automorphism in 
the level $5$ and $6$ cases. The pullback by a generic automorphism
is not isomorphic to the pullback by an exotic isomorphism, since then 
$Y$ itself would be isomorphic to the pullback by an exotic isomorphism.
In the case of level $4$ the quotient is isomorphic to 
the pullback of $Y''$ by a generic automorphism, where $Y''$ has
either level $4$ or $8$, depending upon the choice of subgroup.
From the Kodaira types of the singular fibers neither is isomorphic 
to the pullback of $Y$ by an exotic automorphism. 

It remains to rule out the possibility that a surface in Table 4.1
acquires extra level structure in characteristic $p$ in the form of 
a cyclic \'etale subgroup. Assume first that this group has order 
prime to $p$. There is a forget the extra level structure map between
modular curves for the two moduli problems. It has positive degree
and does not admit a section even when restricted to characteristic 
$p$ which is prime to both level structures. 

Finally we claim that $Y_{\eta }$ has no \'etale subgroup of order $p$
equal to $char(k)$.
If such were to exist, write $E$ for the quotient elliptic curve. Now 
the dual isogeny, $E\to Y_{\eta }$, is purely inseparable and  
may be identified with the relative Frobenius morphism, $F_{E/\eta }$.
This identifies $Y_{\eta }$ with $E^{(p)}$ \cite {Sil, II.2.10}. It 
follows that the $J$-invariant of $Y$ is a $p$-th power. This contradicts
the fact that the singular fibers of $Y$ have Kodaira types I$_r$ with 
$p\nmid r$.
\hfill $\square $\enddemo

\bigpagebreak

\subhead 7. Classification of desingularized fiber products with $h^3=0$\endsubhead 

Let $W$ be  a desingularized fiber product of semi-stable 
elliptic surfaces with vanishing third Betti number. By Corollary 3.2 and \S 4 
$W$ is isomorphic to 
the desingularization of a fiber product, $Y\times _XY'$, constructed as follows: 
$\pi :Y\to X=\Bbb P^1_k$ is a BLIS fibration and $\pi ':Y'\to X$ is the base change of 
a BLIS fibration, $\pi '':Y''\to X$, by an element $\alpha \in Aut(X)$ with 
$\alpha (S)=S''$. (We write $\pi '=\alpha ^*\pi ''$ and $Y'=\alpha ^*Y''$.) 
Furthermore $Y$ and $Y'$ must not be isogenous. The results of \S 5 and \S 6
permit us to list the possibilities. This is done in Table 7.1 below.

\medskip
\proclaim{\hfil Table 7.1. Fiber products of rational elliptic surfaces with $S=S'$}
$$\vbox
{\offinterlineskip\tabskip=0pt\hrule
\halign{ &\vrule#&\strut $\quad #\hfil\quad$\cr
\hr
&\text{characteristic} && \text{Types} &&\text{Kodaira types}  &\cr
\hr
\noalign{\hrule}
\hr
& 2 && 3,5 && \text{I}_3-\text{I}_5, \text{I}_3-\text{I}_5, \text{I}_3-\text{I}_1, \text{I}_3-\text{I}_1   &\cr
\hr
\noalign{\hrule}
\hr
& 2 && 5,5 && \text{I}_5-\text{I}_5, \text{I}_5-\text{I}_1, \text{I}_1-\text{I}_5, \text{I}_1-\text{I}_1   &\cr
\hr
\noalign{\hrule}
\hr
& 3 && 2,2 && \text{I}_4-\text{I}_4, \text{I}_4-\text{I}_2, \text{I}_2-\text{I}_4, \text{I}_2-\text{I}_2 &\cr
\hr
\noalign{\hrule}
\hr
& 5 && 2,6 &&  \text{I}_4-\text{I}_6, \text{I}_2-\text{I}_3, \text{I}_4-\text{I}_2, \text{I}_2-\text{I}_1 &\cr 
\hr
\noalign{\hrule}
\hr
& 5 && 6,6 &&\text{I}_6-\text{I}_2, \text{I}_3-\text{I}_3, \text{I}_2-\text{I}_6, \text{I}_1-\text{I}_1 &\cr 
\hr
\noalign{\hrule}
\hr
& 7 && 2,6 && \text{I}_4-\text{I}_6, \text{I}_4-\text{I}_3, \text{I}_2-\text{I}_2, \text{I}_2-\text{I}_1 &\cr 
\hr
\noalign{\hrule}
\hr
&7 && 6,6 &&  \text{I}_6-\text{I}_3, \text{I}_3-\text{I}_6, \text{I}_2-\text{I}_2,
\text{I}_1-\text{I}_1 &\cr 
\hr
\noalign{\hrule}
\hr
& 11 && 2,5 && \text{I}_4-\text{I}_5, \text{I}_2-\text{I}_5, \text{I}_4-\text{I}_1, \text{I}_2-\text{I}_1 &\cr 
\hr
\noalign{\hrule}
\hr
&11 && 5,5 && \text{I}_5-\text{I}_5, \text{I}_5-\text{I}_5, \text{I}_1-\text{I}_1, \text{I}_1-\text{I}_1 &\cr 
\hr
\noalign{\hrule}
\hr
& 17 && 2,6 && \text{I}_2-\text{I}_6, \text{I}_4-\text{I}_3, \text{I}_4-\text{I}_2, \text{I}_2-\text{I}_1 &\cr
\hr
\noalign{\hrule}
\hr
&17 && 6,6 &&  \text{I}_6-\text{I}_6, \text{I}_3-\text{I}_2, \text{I}_2-\text{I}_3,
\text{I}_1-\text{I}_1 &\cr 
\hr
\noalign{\hrule}
\hr
&29 && 5,6 && \text{ I}_5-\text{I}_6, \text{I}_5-\text{I}_2, \text{I}_1-\text{I}_3, \text{I}_1-\text{I}_1 &\cr
\hr
\noalign{\hrule}
\hr
&31 && 3,5 &&  \text{I}_3-\text{I}_5, \text{I}_3-\text{I}_5, \text{I}_3-\text{I}_1, \text{I}_3-\text{I}_1   &\cr
\hr
\noalign{\hrule}
\hr
&31 && 5,5 && \text{I}_5-\text{I}_5, \text{I}_5-\text{I}_1, \text{I}_1-\text{I}_5, \text{I}_1-\text{I}_1 &\cr 
\hr
\noalign{\hrule}
\hr
&41 && 5,6 && \text{ I}_5-\text{I}_6, \text{I}_5-\text{I}_2, \text{I}_1-\text{I}_3, \text{I}_1-\text{I}_1&\cr
\hr
\noalign{\hrule}
\hr
&73 && 3,6 && \text{I}_3-\text{I}_6, \text{I}_3-\text{I}_3, \text{I}_3-\text{I}_2,
\text{I}_3-\text{I}_1 &\cr
\hr
\noalign{\hrule}
\hr
&73 && 6,6 && \text{I}_6-\text{I}_2, \text{I}_3-\text{I}_6, \text{I}_2-\text{I}_3,
\text{I}_1-\text{I}_1 &\cr 
\hr
\noalign{\hrule}
\hr
&251&& 2,5 &&  \text{I}_2-\text{I}_5, \text{I}_4-\text{I}_5, \text{I}_4-\text{I}_1, \text{I}_2-\text{I}_1&\cr
\hr
\noalign{\hrule}
\hr
&251 && 5,5 && \text{ I}_5-\text{I}_5, \text{I}_5-\text{I}_1, \text{I}_1-\text{I}_5, \text{I}_1-\text{I}_1 &\cr 
\hr
\noalign{\hrule}
\hr
&919&& 5,6 && \text{ I}_5-\text{I}_6, \text{I}_5-\text{I}_3, \text{I}_1-\text{I}_2, \text{I}_1-\text{I}_1 &\cr
\hr
\noalign{\hrule}
\hr
&9001&& 5,6 &&\text{ I}_5-\text{I}_2, \text{I}_5-\text{I}_3, \text{I}_1-\text{I}_6, \text{I}_1-\text{I}_1 &\cr 
\hr
\noalign{\hrule}
\hr
}\hrule}$$
\endproclaim

\example {Explanation of Table 7.1} The rows in the table are in bijective 
correspondence with {\it equivalence classes} of desingularized fiber products
of semi-stable elliptic surfaces with vanishing third Betti number defined
over the algebraic closure of the prime field. Two 
desingularized fiber products, $W_1$ and $W_2$, are said to be {\it equivalent}
if one of the following two relations exists among the pairs of 
semi-stable elliptic surfaces,  $(\pi _1,\pi '_1)$
and $(\pi _2,\pi '_2)$, from which they are constructed: There 
exists $\alpha \in Aut(X)$ such that

\qquad \qquad (i) $\pi _1$ is isogenous to $\alpha ^*\pi _2$  and 
$\pi '_1$ is isogenous to $\alpha ^*\pi '_2$, or

\qquad \qquad (ii) $\pi _1$ is isogenous to $\alpha ^*\pi '_2$ and
$\pi '_1$ is isogenous to $\alpha ^*\pi _2$.

For each equivalence class of desingularized fiber products with
vanishing third Betti number, 
the first two columns of the table indicate the characteristic of the base
field and the types of the BLIS fibrations involved.  The third 
column indicates the nature of the four singular fibers in the fiber 
product. For simplicity this information is only given when the 
exponents of the BLIS fibrations are both $0$ and the level 
is minimal for the given type. A complete list of singular fibers
may be obtained by replacing the $4$-tuple of Kodaira types in the table with
the $4$-tuple of Kodaira types of isogenous elliptic surfaces. For example
applying this to the final line in the table, $p=9001$, gives the 
complete list of singular fibers, $(n,n')\in (\Bbb Z_{\geq 0})^2$:
$$\align  &\text{ I}_{p^n5}-\text{I}_{p^{n'}2},\ \text{I}_{p^n5}-\text{I}_{p^{n'}3},\  \text{I}_{p^n}-\text{I}_{p^{n'}6},\ \text{I}_{p^n}-\text{I}_ {p^{n'}} ,\qquad \text{and}\\
&\text{ I}_{p^n5}-\text{I}_{p^{n'}6},\ \text{I}_{p^n5}-\text{I}_{p^{n'}},\ \text{I}_{p^n}-\text{I}_{p^{n'}2},\ \text{I}_{p^n}-\text{I}_{p^{n'}3}.\endalign $$
\endexample
\example {Justification of Table 7.1} Consider first the case when the 
BLIS fibrations $\pi $ and $\pi ''$ have the same type $\tau $. Then $S=S''$ and 
$\alpha \in Aut(X,S)=Aut(\tau ,p)$, where $p=\op{char}(k)$. By Proposition 6.1
$\alpha $ cannot be generic, but any exceptional $\alpha $ gives rise a desingularized
fiber product $W$ with $H^3(W,\ql )=0$. If $\alpha '\in Aut(\tau ,p)$ is also 
exceptional and $\pi '''$ is the base change of $\pi ''$ with respect to $\alpha '$,
then $\pi $ and $\pi '''$ are isogenous if and only if $Y'$ is isogenous to its 
base change with respect to $\alpha \inv \circ \alpha '$ which is equivalent to 
$\alpha $ and $\alpha '$ lying in the same left coset of $Aut(\tau ,0)$ in 
$Aut(\tau ,p)$. A list of all possible primes $p$ and left cosets of  $Aut(\tau ,0)$ in 
$Aut(\tau ,p)$ is given in Table 5.2. If $Aut(\tau ,0)$ has index two in $Aut(\tau ,p)$,
there is only one non-trivial coset, so there is only one equivalence class of 
desingularized fiber product with $h^3=0$ for the given choice of characteristic
and type. When $Aut(\tau ,0)$ has index three the elements of $Aut(\tau ,p)$ which appear in Table 5.2 have the form $\alpha $ and $\alpha ^2$, where
$\alpha $ has order $3$.  
Base changing the fiber product $Y\times _X\alpha ^*Y$ by $\alpha \inv $ gives
the fiber product $(\alpha ^2)^*Y \times _XY $ which is equivalent to 
$Y \times _X(\alpha ^2)^*Y$. Again there is 
only one equivalence class for the given characteristic and type. 

Suppose now that the types $\tau $ and $\tau ''$ of the BLIS fibrations $\pi $ 
and $\pi ''$ are different. There are only finitely many characteristics $p$ 
in which  
there is an element $\alpha \in Aut(X)$ with $\alpha (S)=S''$. In fact this 
occurs exactly when the cross ratio orbits of $S$ and $S''$ coincide in  
characteristic $p$ and $p\nmid \tau \tau ''$. When $(\tau ,\tau '')=(2,3)$,
these conditions are never met. When $\tau \in \{ 2,3\} $ and 
$\tau "\in \{ 5,6\} $,
then we are dealing with the exceptional primes for types $5$ and $6$ listed in
Table 5.1. To find the characteristics $p$ in the case $\tau =5, \tau '=6$ we refer to Table 4.1 and ask when $\frak K=\frak C$, ie. when 
$\kappa =(11+5\sqrt 5)/(11-5\sqrt 5)\in \{ -8,-\frac 18, 9,\frac 19,\frac 98,\frac 89\} $. It is easy to analyze the various possibilities.
For example if $\kappa =8/9$, then $9(11+5\sqrt 5)=8(11-5\sqrt 5)$, or 
equivalently $11=-85\sqrt 5$. Squaring both sides gives $0=36004=2^2\cdot 9001$.
The case $p=2$ is ruled out, since $2$ divides the level $6$, so $p=9001$ is the
only solution. The case $\kappa =9/8$ also gives $p=9001$, while $\kappa =-8$
or $-1/8$ yields $p=919$ and $\kappa =9$ yields $p=29$ or $41$ as does
$\kappa =1/9$. 

When $\tau <\tau ''$, $Y$ and $Y'$ are not isogenous. This follows from the fact
that if $Y_{\eta }$ has a subgroup of prime order $l$, then so does any isogenous
elliptic curve. To see that each line in Table 7.1 with $\tau <\tau ''$ 
corresponds to a most one equivalence class observe  
that the choice of $\alpha \in Aut(X)$ identifying $S$ with $S''$ is unique 
up to precomposing with elements of $Aut(S,\tau)$. Since the characteristics
which appear in the lines with $\tau <\tau ''$ are never exceptional
for $\tau $ in the sense of Table 5.1, the effect of altering the choice of 
$\alpha $ is to replace $Y$ by an isogenous elliptic surface.
\endexample

\bigpagebreak
\subhead 8. The canonical sheaf and some Hodge numbers \endsubhead

We continue to use the notation of the previous section. In particular $W$ is a 
desingularized fiber product with $H^3(W,\ql )=0$ constructed from
BLIS fibrations $\pi $ and $\pi ''$ in characteristic $p$ of exponents 
$n, n'\in \Bbb Z_{\geq 0}$. Furthermore $\pi '=\alpha ^*\pi ''$ where 
$\alpha \in Aut(X)$ with $\alpha (S)=S''$.

\proclaim {Proposition 8.1} $h^0(W,\Omega ^3_{W/k})=p^n+p^{n'}-1$.
\endproclaim
\demo{Proof} The morphisms $\pi :Y\to X$, $\pi ':Y\to X$, and $\bar f:\bar W\to X$
are projective, local complete intersection morphisms. It follows from 
\cite {Kl, Example 7(ii) and Corollary 19} that relative dualizing 
sheaves $\omega _{Y/X}$, $\omega _{Y'/X}$, and $\omega _{\bar W/X}$ exist, 
that each is an invertible sheaf, and that the restriction to any open 
subscheme which is smooth over $X$ is isomorphic to the determinant of 
the relative K\"ahler differentials. Write $q:\bar W\to Y$ and $q':\bar W\to Y'$
for the projections. From the behaviour of relative K\"ahler 
differentials in fiber products and the fact that the non-smooth locus
is in codimension $2$ we conclude that 
$\omega _{\bar W/X}\simeq q^*\omega _{Y/X}\otimes (q')^*\omega _{Y'/X}$.
A similar argument shows that there is a relative dualizing sheaf, $\omega _{\bar W/k}$,
which is isomorphic to $\bar f^*\Omega _{X/k}\otimes \omega _{\bar W/X}$.
The adjunction formula gives $\Omega ^3_{W/k}\simeq \sigma ^*\omega _{\bar W/k}\otimes \Cal O_W(Q)$, where $Q\subset W$ is the exceptional divisor of the blow-up, $\sigma :W\to \bar W$.
The tautological map, $\pi ^*\pi _*\omega _{Y/X}\to \omega _{Y/X}$, is 
an isomorphism, as can be verified by a local calculation \cite  {Liu, 9.4.3.5}.
Thus 
$$\Omega ^3_{W/k}\simeq f^*(\Omega _{X/k}\otimes \pi _*\omega _{Y/X}\otimes 
\pi '_*\omega _{Y'/X})\otimes \Cal O_W(Q).\tag 8.1$$
Relative duality gives an isomorphism, $\pi _*\omega _{Y/X}^{\vee }\simeq 
R^1\pi _*\Cal O_Y$. Now 
$$12deg(\pi _*\omega _{Y/X})=-12deg(R^1\pi _*\Cal O_Y)=12\left( \chi (\Cal O_X)-
\chi (R^1\pi _*\Cal O_Y)\right) =12\chi (\Cal O_Y)=e(Y),$$
where $\chi (\Cal F)$ is the Euler characteristic of the coherent sheaf $\Cal F$,
$e(Y)$ is the $l$-adic Euler characteristic for $Y$, 
the second to last equality comes from the Leray spectral sequence for $\pi $,
the last equality is Noether's Theorem and the fact that $K_Y^2=0$, which follows from 
$\Omega _{Y/k}^2\simeq \pi ^*(\Omega _{X/k}\otimes \pi _*\omega _{Y/X})$.
Since $\pi $ is 
semi-stable and in particular tamely ramified, $e(Y)$ is the sum of the Euler
characteristics of the singular fibers \cite {Sch,3.2}. Thus $e(Y)=12p^n$.
Since $X\simeq \Bbb P^1_k$, $\pi _*\omega _{Y/X}\simeq \Cal O_{\Bbb P^1}(p^n)$. 
Apply the projection formula and the fact that $f_*\Cal O_W(Q)\simeq \Cal O_{\Bbb P^1}$
to (8.1) to conclude
$$H^0(W,\Omega ^3_{W/k})\simeq H^0(X,f_*\Omega ^3_{W/k})\simeq H^0(\Bbb P^1,\Omega _{\Bbb P^1/k}(p^n+p^{n'}))\simeq k^{p^n+p^{n'}-1}.\qquad \qquad \hfill \square $$ 
\enddemo

\proclaim {Proposition 8.2} $h^1(W,\Cal O_W)=0$, $h^2(W,\Cal O_W)=p^n+p^{n'}-2$, 
$h^3(W,\Cal O_W)=p^n+p^{n'}-1$.
\endproclaim
\demo{Proof} The third assertion follows from the previous proposition and 
Serre duality. 

The first assertion follows from the Leray spectral sequence:
Since $g_X=0$, $H^1(X,f_*\Cal O_W)\simeq H^1(X,\Cal O_X)=0$ and 
$$R^1f_*\Cal O_W\simeq R^1\bar f_*\Cal O_{\bar W}\simeq R^1\pi _*\Cal O_Y\oplus 
R^1\pi '_*\Cal O_{Y'},$$
where the degrees of the two invertible sheaves on the right is negative by the 
proof of the previous proposition. 

For the second assertion observe that $R^2f_*\Cal O_W$ is an invertible sheaf of 
negative degree and hence has no non-zero global sections. Indeed
$$R^2f_*\Cal O_W\simeq R^2\bar f_*\Cal O_W\simeq R^1\pi _*R^1q_*(q')^*\Cal O_{Y'}
\simeq R^1\pi _*\pi ^*R^1\pi '_*\Cal O_{Y'}\simeq R^1\pi _*\Cal O_Y\otimes R^1\pi '_*\Cal O_{Y'},$$
where the second isomorphism comes from the Leray spectral sequence for $\bar f=\pi \circ q$, the third isomorphism is flat base change \cite {Ha, III.9.3} and the fourth
is the projection formula \cite {Ha, III.Ex8.3}. Thus
$$H^2(W,\Cal O_W)\simeq H^1(X,R^1\pi _*\Cal O_Y\oplus R^1\pi '_*\Cal O_{Y'})\simeq 
k^{p^n+p^{n'}-2},$$
since $deg(R^1\pi _*\Cal O_Y)=-e(Y)/12=-p^n$ and $deg(R^1\pi '_*\Cal O_Y)=-p^{n'}$.
\hfill $\square $\enddemo

\bigpagebreak
\subhead 9. Proof of Theorem 1.1\endsubhead 

Part (i) of the theorem was established in the Justification of Table 7.1.

For part (ii) recall from \S 7 that $W$ is constructed from data 
$(\pi , \pi '', \alpha )$, where $\pi $ and $\pi ''$ are BLIS fibrations
and $\alpha \in Aut(X)$ satisfies $\alpha (S)=S''$. 
Note that every BLIS fibration is defined over the prime field.
If $\pi $ and $\pi ''$ have the same type, then $\alpha $ appears
in Table 5.2. It is defined over the prime field except when $char(k)=2$,
in which case it is defined over the field with $4$ elements.
Suppose the types satisfy $\tau <\tau ''$. In the case of the first line
of Table 7.1 we may take $\alpha =Id$. In the remaining cases $\alpha $
is defined over the prime field since the cross ratio of $S$ is. Finally 
the blow-up of 
the fiber product along $\bar W_{sing}$, $\sigma :W\to \bar W$, is 
defined over the field of definition of $\bar W$. 

Part (iii) follows from Proposition 8.1.

Part (iv) follows from Proposition 8.1 and the finiteness of the number of 
lines in Table 7.1.

\bigpagebreak
\subhead 10. Projective threefolds with trivial canonical sheaf and $h^3=0$\endsubhead

It is remarkable that there exist non-singular projective threefolds, $V$,  
with trivial canonical sheaf and $H^3(V,\ql )=0$. Such varieties do not exist 
in characteristic zero by Hodge theory. The first example was constructed 
by Hirokado \cite {Hi1} in characteristic $3$. More recently 
Schr\"oer \cite  {Schr} and then Ekedahl \cite {Ek} have constructed
a few further examples in characteristics $2$ and $3$; some of these have 
moduli. 

\proclaim {Proposition 10.1}
(i) The fiber product of level $(4,4)$ in characteristic $3$ 
(see Table 7.1 line 3)
admits several non-isomorphic small resolutions which are projective threefolds 
with trivial canonical sheaf and $h^3=0$.

(ii) The fiber products of other types of rational elliptic surfaces
in Table 7.1 do not admit small projective resolutions. Small resolutions
do exist as algebraic spaces.

(iii) There exist projective threefolds with trivial canonical sheaf and $h^3=0$
which admit a pencil of Kummer surfaces in which the general fiber is not 
supersingular. 
\endproclaim 

\demo{Proof}
(i) In characteristic $3$ consider the fiber product of the level $4$ surface in Table 4.1
with its base change by an exceptional automorphism (Table 5.1 line 1 and 
Table 7.1 line 3). Every irreducible component of every fiber of the tautological
map, $\bar f:\bar W\to X$, is non-singular. Ordering the components of the
singular fibers and then successively blowing up one component at 
a time in the order chosen gives rise to a morphism, $\gamma :\widehat W\to \bar W$,
which collapses finitely many rational curves in $\widehat W$. The variety
$\widehat W$ is projective, non-singular, and has trivial canonical sheaf.
There is a birational morphism, $\varpi :W\to \widehat W$, which on every component,
$\Bbb P^1\times \Bbb P^1$, of the exceptional divisor $Q\subset W$ is projection
onto one factor. It follows from the Leray spectral sequence for $\varpi $ that
$h^3(\widehat W,\ql )=0$. 

To verify that the isomorphism class of $\widehat W$ depends on the 
order in which the blow-ups are performed, observe first that each singular 
fiber of $\bar f$  is the product of a Kodaira type I$_2$ fiber with a 
Kodaira type I$_4$ fiber or is the self-product of an I$_2$ fiber or an
I$_4$ fiber. The components of I$_2$ and I$_4$ which meet the identity
section will be called identity components. These components together
with the component of I$_4$ which does not meet the identity component
will be called {\it even}. For $s\in S$ a component of $\bar f\inv (s)$
is called even if it is the product of even components of $\pi \inv (s)$
and $(\pi ')\inv (s)$. Every singular point
of $\bar W$ is contained in an even component of a singular fiber. Thus
blowing up the even components already gives a non-singular variety, $\widehat W_{even}$. 
Blowing up the strict transforms of the remaining components has no 
effect, as they are Cartier divisors. Since the even components 
are disjoint, the isomorphism class of $\widehat W_{even}$ is independent
of the order in which the even components are blown up. The inversion maps 
on the generic fibers of $\pi $ and $\pi '$ extend to  biregular
involutions of $Y$ and $Y'$. The product gives a biregular involution of $\bar W$
which lifts to a biregular involution, $\iota $, of $\widehat W_{even}$. To construct
a resolution $\widehat W\to \bar W$ such that inversion is not biregular fix a 
singular fiber of $\bar W$ of type I$_2-$I$_4$. Blow-up the component
which is the product of the identity component of the I$_2$ fiber with
a component of the I$_4$ fiber which does not meet any section of order $2$.
Then blow up the component which meets the identity section of $\bar f$.
No matter in what order the remaining fiber components are blown up,
inversion will not give a regular involution of the resulting variety,
$\widehat W$. 

(ii) The obstruction is the existence of a fiber of Kodaira type $I_1$
in either $Y$ or $Y'$. The proof of \cite {Sch2, 3.1(iii)} goes through in this case.
The assumption in {\it loc. cit.} that $S\neq S'$ does not hold here,
but all one needs is that $Y$ and $Y'$ are not isogenous. There is no
obstruction to small resolution in the category of algebraic spaces.

(iii) Write $\tilde W_{even}$ for the blow-up of $\widehat W_{even}$ along the 
fixed point locus of $\iota $. Then $\iota $ lifts to a biregular 
involution $\tilde {\iota }\in Aut(\tilde W_{even})$. Denote the 
quotient by $\frak W$. Since the fixed locus of $\tilde {\iota }$ is a non-singular divisor
and the characteristic is different from $2$, $\frak W$ is non-singular.
Any non-zero global section of $\Omega ^3_{\tilde W_{even}/k}$ descends to give a 
nowhere vanishing section of $\Omega ^3_{\frak W/k}$. Blowing up $\widehat W_{even}$
along $16$ disjoint $\Bbb P^1$'s to produce $\tilde W_{even}$ does not change the 
third Betti number. The map $H^3(\frak W,\ql )\to H^3(\tilde W_{even},\ql )$ is 
injective as one can see by applying the projection formula \cite {Mi, Vi.11.6a,d} with one factor 
being $1\in H^0(\tilde W_{even},\ql )$. Thus $H^3(\frak W,\ql )=0$.
\hfill $\square $\enddemo

\example {Open question 10.2}
(i) Do there exist non-singular projective threefolds with trivial canonical 
sheaf and $h^3=0$ in characteristics other than $2$ and $3$?
\endexample

\bigpagebreak
\subhead 11. Deformations and lifting to characteristic zero \endsubhead 

The smooth projective threefolds with trivial canonical sheaf constructed by 
Hirokado, Schr\"oer, and Ekedahl all have obstructed deformations.
This manifests itself when one tries to lift these varieties to 
characteristic zero \cite {Hi1}, \cite {Schr}, \cite {Ek}. This result 
stands in notable contrast to the situation in characteristic zero where
the deformation theory of smooth, projective threefolds with trivial
canonical sheaf is unobstructed.

\proclaim {Proposition 11.1}
Let $W$ be a desingularized fiber product of semi-stable rational elliptic 
surfaces with $h^3(W,\ql )=0$. Then $W$ does not admit a formal lifting
to characteristic zero.
\endproclaim
\demo{Proof} The argument in \cite {Schr, \S 2} works here with minor
modification. By Proposition 8.2 $H^1(W,\Cal O_W)=0$. The assumption
that the elliptic surfaces are rational (ie. exponent $0$)  implies
$H^2(W,\Cal O_W)=0$. As in {\it loc. cit.} it follows that if a formal
lifting were to exist, then such a lifting would exist over a
Noetherian, integral base, $B$, and it would be isomorphic to a formal 
completion of a smooth projective scheme, $\bold f:\Cal W\to B$. 
Furthermore restriction would give an isomorphism, $Pic(\Cal W)\to Pic(W)$.
Since $\Omega ^3_{W/k}\simeq \Cal O_W(Q)$ and $H^i(Q,\Cal N_{Q/W})=0$
for all $i$, the relative Hilbert scheme, $Hilb_{\Cal W/B}$, is 
locally \'etale over $B$ at $Q$. Thus $Q$ lifts to a unique effective divisor 
$\Cal Q\subset \Cal W$. The isomorphism of Picard groups gives that
$\Omega _{\Cal W/B}^3\simeq \Cal O_{\Cal W}(\Cal Q)$. Passing to the
fraction field $K$ of $B$ gives a smooth, projective threefold, $\Cal W_K$ over a
field of characteristic zero with $h^0(\Cal W_K, \Omega ^3_{\Cal W_K/K})\neq 0$. By base change theorems, $R^3\bold f_*\Bbb Z/l^n$, is a constant
sheaf on $B$ for all $n$. Thus $h^3(\Cal W_{\bar K},\ql )= 0$. Passing to 
an associated complex manifold and applying Hodge theory gives a contradiction.
\hfill $\square $\enddemo

The next result shows that the threefolds, $\widehat W$, constructed in the previous section
are rigid. Define the tangent sheaf, $\Cal T_{\widehat W/k}:=\Cal Hom _{\Cal O_{\widehat W}}(\Omega _{\widehat W/k},\Cal O_{\widehat W})$. 
\proclaim {Theorem 11.2} Let $\widehat W$ be a small projective resolution
of a fiber product of semi-stable elliptic surfaces with section and at least
one singular fiber. Suppose that $\Omega ^3_{\widehat W/k}\simeq \Cal O_{\widehat W}$. If $h^3(\widehat W,\ql )=0$, then $H^1(\widehat W,\Cal T_{\widehat W/k})=0$.
\endproclaim
\demo{Proof} By the triviality of the canonical sheaf there is an isomorphism,
$\Omega ^2_{\widehat W/k}\simeq \Cal T_{\widehat W/k}$. By Serre duality it
suffices to show $H^2(\widehat W,\Omega ^1_{\widehat W/k})=0$. As verifying
this seems no easier than computing all the Hodge numbers of $\widehat W$ 
we prove
\proclaim {Proposition 11.3} The only non-zero Hodge numbers of $\widehat W$ are 
$h^{00}(\widehat W)=h^{33}(\widehat W)=1$, 
$h^{03}(\widehat W)=h^{30}(\widehat W)=1$, and $h^{11}(\widehat W)=
h^{22}(\widehat W)=36$.
\endproclaim
\demo{Proof} By Proposition 8.2 and Serre duality it suffices to compute
$H^i(\widehat W,\Omega _{\widehat W/k})$ for each $i$. As in \cite {Sch3, \S 11} the cohomology of $\Omega _{\widehat W/k}$ may be computed using the Leray spectral sequence for $\hat f:\widehat W\to X$. In particular $H^0(\widehat W,\Omega _{\widehat W/k})=0$, since the natural map, $\Omega _X\to \hat f_*\Omega _{\widehat W/k}$ is an isomorphism \cite {Sch3, 11.8(iv)}. Write $m_s$ (respectively $m'_s$) for the number of irreducible components in $\pi \inv (s)$ (respectively $(\pi ')\inv (s)$). The main technical lemma we use is: 
\proclaim {Lemma 11.4} $R^1\hat f_*\Omega _{\widehat W/X}\simeq 
\Cal O_X^3\oplus \Cal O_X(-S) \oplus _{s\in S}i_{s*}\Cal O_s^{m_sm'_s-1}$.
\endproclaim
The proof of the lemma will be given at the end of the section. To use it to 
prove the proposition, consider the short exact sequence 
$$0\to \hat f^*\Omega _{X/k}\to \Omega _{\widehat W/k}\to \Omega _{\widehat W/X}\to 0$$
and the corresponding long exact sequence of derived functors,
$$R^1\hat f_*\Omega _{\widehat W/k}@>c>>R^1\hat f_*\Omega _{\widehat W/X}
@>e>> R^2\hat f_*\hat f^*\Omega _{X/k}\to  R^2\hat f_*\Omega _{\widehat W/k}
\to R^2\hat f_*\Omega _{\widehat W/X}.\tag 11.5$$
Write $\Cal L$ for $R^1\pi _*\Cal O_Y$ and $\Cal L'$ for $R^1\pi '_*\Cal O_{Y'}$.
By \cite {Sch3, 11.8(iii)} and the fact that $Y$ and $Y'$ are rational elliptic surfaces,
$$R^2\hat f_*\hat f^*\Omega _{X/k}\simeq \Omega _{X/k}\otimes R^2\hat f_*\Cal O _{\widehat W}\simeq \Omega _{X/k}\otimes \Cal L\otimes \Cal L'\simeq \Cal O_{\Bbb P^1}(-4). $$
Furthermore $c$ is injective since 
$$0\to \Omega _{X/k}@>\alpha >> \hat f_*\Omega _{\widehat W/k}\to 
\hat f_*\Omega _{\widehat W/X}@>e''>>
\Omega _X\otimes R^1\hat f_*\Cal O_{\widehat W}\to \op{Ker}(c)\to 0$$
is exact, $\alpha $ is an isomorphism, and  
$$\hat f_*\Omega _{\widehat W/X}\simeq \Cal L\inv (-S)\oplus (\Cal L')\inv (-S)
\simeq \Cal O_{\Bbb P^1}(-3)^{\oplus 2}\simeq \Omega _X\otimes (\Cal L\oplus \Cal L')\simeq \Omega _X\otimes R^1\hat f_*\Cal O_{\widehat W}$$
by \cite {Sch3, 11.8(ii) and 11.14(i)}.

As the ranks of the first two $\Cal O_X$-modules in (11.5) are $3$ and $4$,
the map $e$ is non-zero. By Lemma 11.4 it must be the projection of $R^1\hat f_*\Omega _{\widehat W/X}$ on the direct summand $\Cal O_X(-S)\simeq \Cal O_{\Bbb P^1}(-4)$.
Thus 
$$R^1\hat f_*\Omega _{\widehat W/k}\simeq \Cal O_X^{\oplus 3}
\oplus _{s\in S}i_{s*}\Cal O_s^{m_sm'_s-1},$$
where $i_s:s\to X$ is the inclusion. Consequently,  
$$h^1(\widehat W, \Omega _{\widehat W/k})=h^1(X,\hat f_*\Omega _{\widehat W/k})
+ h^0(X,R^1\hat f_*\Omega _{\widehat W/k})=1+3+\sum _{s\in S}(m_sm'_s-1).$$
By Proposition 10.1 and line 3 of Table 7.1 $\sum _{s\in S}m_sm'_s=16+8+8+4=36$. Thus $h^1(\widehat W, \Omega _{\widehat W/k})=36$. Now (11.5) and \cite {Sch3, 11.15} yield 
$$ R^2\hat f_*\Omega _{\widehat W/k}
\simeq R^2\hat f_*\Omega _{\widehat W/X}\simeq \Cal L\oplus \Cal L'
\simeq \Cal O_{\Bbb P^1}(-1)^{\oplus 2}.$$
Thus $h^3(\widehat W,\Omega _{\widehat W/k})=h^1(X,R^2\hat f_*\Omega _{\widehat W/k})=0$ and 
$$h^2(\widehat W,\Omega _{\widehat W/k})=h^0(X,R^2\hat f_*\Omega _{\widehat W/k})
+h^1(X,R^1\hat f_*\Omega _{\widehat W/k})=0.$$

It remains only to prove the lemma.

\demo{Proof of Lemma 11.4}
Observe that the sheaf $\Omega _{\widehat W/X}$ is torsion free and hence flat over $X$. This follows from the fact that for a closed point $w\in \widehat W$ with image $x\in X$ there are local parameters $x_1,...,x_3\in \Cal O_{\widehat W,w}$ such that the pullback of a uniformizing parameter $t\in \Cal O_{X,x}$ satisfies $t=x_1...x_i$ for some $i\leq 3$ which implies that $dt\in \Omega _{\widehat W/k. w}\simeq \oplus _{i=1}^3\Cal O_{\widehat W,w}dx_i$ is an indivisible element \cite {Ha, II.8.3A}. 

Write $E\subset \widehat W$ for the exceptional locus of the small projective
resolution of singularities, 
$\gamma :\widehat W\to \bar W$. In the exact sequence 
\cite {Sch3, 11.11-12}
$$0\to R^1\hat f_*\gamma ^*\Omega _{\bar W/X}@>\alpha >>
R^1\hat f_*\Omega _{\widehat W/X}\to R^1\hat f_*\Omega _{E/k}\to R^2\hat f_*\gamma ^*\Omega _{\bar W/X}\to R^2\hat f_*\Omega _{\widehat W/X}\to 0,$$
we have
$$\align
R^1\hat f_*\Omega _{E/k}&\simeq \oplus _{s\in S}i_{s*}\Cal O_s^{\oplus m_sm'_s}=:\mu , \\ 
R^1\hat f_*\gamma ^*\Omega _{\bar W/X}&\simeq \Cal O_X^2\oplus \Cal O_X^2(-S)
\oplus \nu \qquad \qquad \qquad 
\text{\cite {Sch3, 11.11-13}} \\
R^2\hat f_*\gamma ^*\Omega _{\bar W/X}&\simeq \Cal L\oplus \Cal L'
\oplus \nu , \qquad \qquad \qquad \qquad \qquad 
\text{\cite {Sch3, 11.13v}}\endalign $$
where $\nu =\oplus _{s\in S}i_{s*}\Cal O_s^{\oplus m_s+m'_s-2}$.  The direct summand $\Cal O_X^2\subset 
R^1\hat f_*\gamma ^*\Omega _{\bar W/X}$ is generated by the Hodge 
cohomology classes of the inverse images of the identity sections of 
$\pi $ and $\pi '$ via the projections $q:\widehat W\to Y$ and 
$q':\widehat W\to Y'$. It is also a direct summand of $R^1\hat f_*\Omega _{\widehat W/X}$. As 
$R^2\hat f_*\Omega _{\widehat W/X}\simeq \Cal L\oplus \Cal L'$ 
\cite {Sch3, 11.15}, the above exact sequence gives rise to an exact 
sequence,
$$0\to \Cal O_X^2(-S)\oplus \nu \to R^1\hat f_*\Omega _{\widehat W/X}/\Cal O_X^2
\to \mu \to \nu \to 0.$$

Set $F=\hat f\inv (s)$. By \cite {Sch3,11.15} $R^2\hat f_*\Omega _{\widehat W/X}$ is locally free and the natural map, $R^2\hat f_*\Omega _{\widehat W/X}\otimes i_{s*}\Cal O_s\to H^2(F, \Omega _{F/k})$, 
is an isomorphism. It follows that the natural map,
$R^1\hat f_*\Omega _{\widehat W/X}\otimes i_{s*}\Cal O_s
\to H^1(F, \Omega _{F/k})$ is also an isomorphism \cite {Ha, 12.11 and 8.2a}. 
For $s\in S$, $dim_k(H^1(F, \Omega _{F/k}))=m_sm'_s+3$ \cite {Sch4}.
Since the rank of $R^1\hat f_*\Omega _{\widehat W/X}$ is four, the torsion
subsheaf, $\tau \subset R^1\hat f_*\Omega _{\widehat W/X}$ satisfies,
$\tau \otimes i_{s*}\Cal O_s\simeq i_{s*}\Cal O_s^{\oplus m_sm'_s-1}$.
The Hodge cohomology classes of the irreducible components  of $F$
generate a subsheaf $\theta _1\subset R^1\hat f_*\Omega _{\widehat W/X}$ isomorphic to $\oplus _{s\in S}i_{s*}\Cal O_s^{\oplus m_sm'_s-1}$. The $-1$
appears in the exponent because 
the class of $F$ is zero. The composition, $\theta _1\to \tau \to 
\tau \otimes i_{s*}\Cal O_s$, is an isomorphism. It follows that 
$\theta _1\simeq \tau $. Thus the previous exact sequence yields the short exact sequence,
$$0\to \Cal O_X^2(-S)@>>> R^1\hat f_*\Omega _{\widehat W/X}/(\Cal O_X^2\oplus \tau )@>>>\oplus _{s\in S}i_{s*}\Cal O_s\to 0.$$
The map $e$ in (11.5) gives rise to the composition,
$$\Cal O_X^2(-S)@>>>R^1\hat f_*\Omega _{\widehat W/X}/(\Cal O_X^2\oplus \tau ) \to \Omega _X\otimes \Cal L\otimes \Cal L',$$
which is non-zero. In the case at hand
$$\Cal O_X(-S)\simeq \Cal O_{\Bbb P^1}(-4)\simeq \Omega _X\otimes \Cal L\otimes \Cal L',$$
so this map is in fact a split surjection. Now the lemma follows from
the short exact sequence,
$$0\to \Cal O_X(-S)@>>> R^1\hat f_*\Omega _{\widehat W/X}/(\Cal O_X(-S)\oplus \Cal O_X^2\oplus \tau )@>>>\oplus _{s\in S}i_{s*}\Cal O_s\to 0.$$
\hfill $\square \square \square $\enddemo\enddemo \enddemo

\example {Remark 11.6} It follows from Proposition 12.1 below that the 
Picard group of $\widehat W$ is free of rank $35=h^{1,1}(\widehat W)-1$.
\endexample 
\bigpagebreak
\subhead 12. Supersingular threefolds \endsubhead

There is a well established notion of supersingular abelian variety in positive 
characteristic. If the dimension is at least two and the base field is algebraically
closed, an abelian variety, $A$, is supersingular if and only if the rank of the 
N\'eron-Severi group, $\rho (A)$, equals $h^2(A,\ql )$. Shioda  \cite {Sh}
has suggested that a smooth projective surface $V$ over an algebraically closed 
field be called supersingular if $\rho (V)=h^2(V,\ql )$. Every unirational surface has this property. Shioda asks if conversely every supersingular surface 
with $\pi _1(V)=\{ 1\} $ is unirational. It seems interesting to consider
these issues in higher dimensions as well.

\proclaim {Proposition 12.1}
Let $W$ be a desingularized fiber product of semi-stable elliptic surfaces with 
vanishing third Betti number. Then

(i) $\pi _1(W)=\{ 1\} $.

(ii) $\rho (W)=h^2(W,\ql )$.
\endproclaim
\demo{Proof} (i) First note that $\pi _1(Y)=\pi _1(Y')=\{ 1\} $, for the elliptic surfaces
from which $W$ is constructed. This is clear if $Y$ is rational (ie. if the exponent 
is $0$). In the general case one need only note that $Y$ contains an open dense subset which is a purely inseparable finite cover of the variety, $U$, which is obtained from a BLIS fibration of exponent $0$ by removing the finite set of points where the fibration is not smooth \cite {Mi, I.5.2h}, \cite {Gr, IX.4.10}. Write $i :\dot X:=X-S=X-S'\to X$ for the inclusion and denote base
change by $i$ by adding a $\dot {}$ to the notation.
For a closed point $x\in \dot X$, there is an exact sequence,
$$\pi _1(f\inv (x))@>>>\pi _1(\dot W)\to \pi _1(\dot X)\to \{ 1\} ,$$
and analogous exact sequences for $\dot Y$ and $\dot Y'$ \cite {Gr, X.1.4} . It follows that 
$\pi _1(\dot W)$ is generated by the images of the groups $\pi _1(\dot Y)$
and $\pi _1(\dot Y')$ under the maps given by the identity sections. 
The surjective map, $\pi _1(\dot W)\to \pi _1(W)$, must be zero, since 
the composition, $\pi _1(\dot Y)\to \pi _1(\dot W)\to \pi _1(W)$, factors
through $\pi _1(Y)$ and the same holds for $Y'$.

(ii) Every desingularized fiber product of semi-stable elliptic surfaces has the property
that its $l$-adic euler characteristic is the sum of the $l$-adic euler characteristics
of the singular fibers \cite {Sch, 8.13}. Write $n$ and $n'$ for the exponents of the 
semi-stable elliptic surfaces $Y$ and $Y'$ used to construct $W$. It is easy to 
compute the $l$-adic euler characteristics, $e(\bar W)=\sum _{s\in S}p^nm_sp^{n'}m'_s$ and 
$e(W)=4\sum _{s\in S}p^nm_sp^{n'}m'_s$. Thus
$$h^2(W,\ql )=\frac 12(e(W)-2+2h^1(W)+h^3(W))=-1+2\sum _{s\in S}m_sp^nm'_sp^{n'}.$$
The rank of the N\'eron-Severi group, which is equal to the Picard group, may
be computed by the localization sequence as the sum of the rank of the Picard
group of the generic fiber plus the contribution from the components of the 
closed fibers \cite {Sch2, 3.2}. The contribution from the generic fiber is
$2$ since $Y$ and $Y'$ are not isogenous and each has Mordell-Weil rank zero
\cite {Be1, Section 4.A}. The only relation of rational equivalence
among the components of the closed fibers is that all fibers are rationally 
equivalent. This fact may be proved by intersecting fiber components 
with the fiber components of a sufficiently general very ample hypersurface
in $W$ and applying a standard result on the intersection matrix of the 
fiber components of a fibered surface \cite {Liu, 9.1.23}.
It follows that the contribution from the closed fibers is 
$2\sum _{s\in S}m_sp^nm'_sp^{n'}-3$. Thus $\rho (W)=h^2(W,\ql )$.
\hfill $\square $\enddemo

The proposition shows that there is no naive obstruction in $l$-adic cohomology to 
the unirationality of $W$. However it is not even known if $W$ is uniruled.
By contrast the threefold of Hirokado is unirational as it was constructed by desingularizing 
a quotient of $\Bbb P^3_{\Bbb F_3}$ by the action of an appropriate vector field \cite {Hi1}. 
The threefolds of Schr\"oer and Ekedahl admit pencils of  
$K3$ surfaces with Picard number $22$, so-called supersingular
$K3$ surfaces, which in characteristic $3$ are Kummer surfaces. These 
are uniruled since supersingular $K3$ surfaces over an algebraically closed
field of characteristic $2$ are unirational \cite {Ru-Sha}  
as are supersingular Kummer surfaces over an algebraically closed field of 
characteristic $3$ \cite  {Sh, Proposition 8}.

Unirationality or uniruledness of $W$ would have consequences for the 
structure of the Chow group. Since $h^3(W,\ql )=0$, the $l$-adic 
Abel-Jacobi map is zero. Is the Chow group of nullhomologous 
one cycles also zero?

\bigpagebreak
\subhead 13. Arithmetic degeneration of rigid Calabi-Yau threefolds\endsubhead 

Let $V$ be a smooth, projective threefold with trivial canonical sheaf and no
first order deformations defined over a field of characteristic zero. We may and will assume that $V$ is defined over a 
number field, $K$. By deformation theory $0=H^1(V,T_{V/K})\simeq H^2(V,\Omega _{V/K})$.
By Hodge theory and comparison theorems, $dim(H^3(V_{\bar K},\ql ))=2$. 
It is interesting to ask how the degenerations of $V$ might relate
to the classical reduction theory of elliptic curves. 
A crude form of this reduction theory may be
stated as follows: Let $\frak o$ be the strict henselisation of the integers
of a number field $K$ at a place $\frak p$ and let $\bold K$ be the fraction 
field of $\frak o$. Let $\bar {\bold K}$ be an algebraic closure of $\bold K$ and
let $I\subset Gal(\bar {\bold K}/\bold K)$ denote the inertia group. Let 
$l$ be a prime distinct from the residue characteristic.

\proclaim {Theorem 13.1} An elliptic curve $E/\bold K$ has a relatively 
minimal regular projective model $\Cal E$ over $\frak o$, which is unique up to isomorphism.

(i) If $I$ acts trivially on $H^{\bullet }(E_{\bar K},\ql )$, then $\Cal E$ is smooth
over $\frak o$.

(ii) If $I$ acts non-trivially, but unipotently, then the closed fiber of 
$\Cal E$ is an $n$-gon of $\Bbb P^1$'s.

(iii) After replacing $\frak o$ by a finite extension if necessary, either 
(i) or (ii) holds.
\endproclaim

A number of rigid threefolds with trivial canonical sheaf may be constructed as 
fiber products of rational semi-stable elliptic surfaces 
\cite {Sch\"u} and \cite {Sch2, \S 7}. The few examples studied so far 
admit a smooth, projective regular model at places
where the inertia group acts trivially even when the elliptic fibrations
degenerate. At places where the inertial action 
on $H^{\bullet }(V_{\bar {\bold K}},\ql )$ is non-trivial, but unipotent, 
one is frequently able to construct a projective regular model in which 
a desingularized fiber product of semi-stable rational elliptic surfaces
with vanishing third Betti number appears as a component of the closed fiber.
A first impression is that such varieties may be playing a role in
three dimensions which is analogous to the role played by 
$n$-gons of $\Bbb P^1$'s in dimension $1$. The author hopes to return to
this question in the future.

\example {Remark 13.2}
It is interesting to note that reducing rigid Calabi-Yau varieties 
(or related varieties) modulo
certain primes of bad reduction is a technique for producing smooth 
projective threefolds over a finite field which do not lift to 
characteristic zero. 
\endexample

\bigpagebreak
\subhead 14. Fiber products of more general elliptic surfaces\endsubhead 

Fiber products of semi-stable elliptic surfaces are attractive because it is very 
easy to resolve the singularities. On the other hand it has become apparent that fiber products of non-semi-stable elliptic surfaces give further insight into the 
issues raised in the previous four sections. For example, such fiber products may be used to construct a one dimensional family of smooth, projective threefolds with trivial canonical sheaf and $h^3=0$.  In the context of \S 12 the construction yields examples of threefolds which can be shown to be inseparably uniruled although not separably uniruled. This construction also arises naturally in the study of degenerations of rigid Calabi-Yau threefolds. The author hopes to say more about these more general fiber products in forthcoming work. Fiber products involving quasi-elliptic surfaces have been studied by Hirokado \cite {Hi2}.

\subhead References \endsubhead 

[Be1] Beauville, A., Le nombre minimum de fibres singuli\`eres d'une courbe stable 
sur $\Bbb P^1$, in S\'eminaire sur les pinceaux de courbes de genre au moins deux,
ed. L. Szpiro, Ast\'erisque 86, p. 97-108 (1981)

[Be2] Beauville, A., Les familles stables de courbes elliptiques sur $\Bbb P^1$
admettant quatre fibres singuli\`eres, C. R. Acad. Sc. Paris, s\'er. I, {\bf 294},
p. 657-660, (1982).

[De] Deligne, P., La conjecture de Weil: II, Pub. math. I.H.\'E.S. {\bf 52},
p. 137-252 (1980)

[Ek] Ekedahl, T., On non-liftable Calabi-Yau threefolds, arXiv:math.AG/0306435 v2 26 May 2004

[Fa-Scha-W\"u] Faltings, G., Schappacher, N., W\"ustholz, G. et al,
Rational Points, $3^{rd}$ enlarged edition, Max-Planck-Institut f\"ur Mathematik, Bonn
(1992)

[Gr] Grothendieck, A., Rev\^etements \'etales et groupe fondamental (SGA I),
Lecture Notes in Mathematics 224, Springer-Verlag, Berlin (1971)

[Ha] Hartshorne, R., Algebraic Geometry, Springer-Verlag, New York (1977)

[Hi1] Hirokado, M., A non-liftable calabi-yau threefold in characteristic 3,
Tohoku Math. J. {\bf 51}, p. 479-487 (1999)

[Hi2] Hirokado, M., Calabi-Yau threefolds obtained as fiber products of
elliptic and quasi-elliptic rational surfaces.  J. Pure Appl. Algebra
{\bf 162}, p. 251-271  (2001)

[I1] Ito, H., On unirationality of extremal elliptic surfaces, Math. Ann. {\bf 310},
717-733 (1998)

[I2] Ito, H., On extremal elliptic surfaces in characteristic $2$ and $3$,
Hiroshima Math. J. {\bf 32}, 179-188 (2002)

[Ka] Katz, N., Moments, Monodromy, and Perversity: A Diophantine Perspective,
Princeton University Press, Princeton (2005)

[Kl] Kleiman, S., Relative duality for quasi-coherent sheaves, Compositio Math. {\bf 41},
39-60 (1980)

[La] Lang, W., Extremal rational elliptic surfaces in characteristic $p$. I:
Beauville surfaces, Math. Z. {\bf 207}, p. 429-437 (1991)

[Liu] Liu, Q., Algebraic Geometry and Arithmetic Curves, Oxford Univeristy Press,
Oxford (2002)

[Mi] Milne, J., \'Etale Cohomology, Princeton University Press (1980)

[Mi-Pe] Miranda, R. and Persson, U., On extremal rational elliptic surfaces,
Math. Z. {\bf 193}, 537-558 (1986)

[Ru-Sha] Rudakov, A. N. and Safarevic. I.R. Supersingular ${K}3$ surfaces over fields of characteristic $2$, Izv. Akad. Nauk SSSR Ser. Mat. 42 (1978), no. 4, 848-869

[Sch] Schoen, C., Torsion in the cohomology of fiber products of elliptic surfaces

[Sch2] Schoen, C., On fiber products of rational elliptic surfaces with section, 
Math. Z. {\bf 197}, 177-199 (1988)

[Sch3] Schoen, C., Complex varieties for which the Chow group mod $n$ is not finite, J. Algebraic Geometry {\bf 11}, 41-100 (2002)

[Sch4] Schoen, C., Invariants of certain normal crossing surfaces, in 
preparation

[Schr] Schr\"oer, S., Some Calabi-Yau threefolds with obstructed deformations over the Witt
vectors,  Compos. Math.  140  (2004), 1579--1592.

[Sch\"u] Sch\"utt, M., New examples of modular rigid Calabi-Yau threefolds, 
Collect. Math.  55  (2004), 219--228.

[Schw] Schweizer, A., Extremal elliptic surfaces in characteristic $2$ and $3$,
manuscripta math. {\bf 102}, 505-521 (2000)

[Sh] Shioda, T., On unirationality of supersingular surfaces, Math. Ann. {\bf 225},
155-159 (1977)

[Sil] Silverman, J., The Arithmetic of Elliptic Curves, Springer-Verlag, New York, (1986).

[Sil2] Silverman, J., Advanced Topics in the Arithmetic of Elliptic Curves,
Springer-Verlag, New York, (1994).
\enddocument